\documentclass{article}

\newtheorem{fed}{\textbf{Definition}}[section]
\newtheorem{thm}[fed]{\textbf{Theorem}}
\newtheorem{lemma}[fed]{\textbf{Lemma}}

\usepackage{amssymb,bbm,graphicx,epsfig,psfrag,epic,eepic,latexsym}
\usepackage{amsmath}
\usepackage{mathrsfs}
\usepackage[dvips]{color}

\begin{document}
\title{V-shaped action functional with delay}
\author{Urs Frauenfelder} \maketitle

\begin{abstract}
In this note we introduce the V-shaped action function with delay in a symplectization
which is an
intermediate action functional between Rabinowitz action functional
and the V-shaped action functional. It lives on the same space as the
V-shaped action functional but its gradient flow equation is a delay equation
as in the case of Rabinowitz action functional. We show that there is a smooth interpolation
between the V-shaped action functional and the V-shaped action functional with delay
during which the critical points and its actions are fixed. On the other hand there
is a bijection between gradient flow lines of the V-shaped action functional with delay
and the ones of Rabinowitz action functional. 
\end{abstract}

\section{Introduction}

V-shaped symplectic homology was introduced in \cite{cieliebak-frauenfelder-oancea},
jointly by the author with Cieliebak and Oancea where it was shown that it is isomorphic to Rabinowitz Floer homology\cite{cieliebak-frauenfelder} and fits
into a long exact sequence with symplectic homology and cohomology. This in particular
allowed to compute Rabinowitz Floer homology of cotangent bundles. It turned out
that these computations coincide with computations in Tate Hochschild cohomology
by Rivera and Wang in \cite{rivera-wang}, where they conjectured a connection
with the algebraic structures discovered in V-shaped symplectic homology by
Cieliebak and Oancea in \cite{cieliebak-oancea}. That there is a close connection
between Tate homology and Rabinowitz Floer homology is in fact to expect in view
of the fact that Rabinowitz action functional is antiinvariant under time reversal.
On the other hand the V-shaped action functional is not antiinvariant anymore under time
reversal and therefore although one has Poincar\'e duality on homology level it gets
lost on chain level. 
\\ \\
The isomorphism between Rabinowitz Floer homology and V-shaped Symplectic homology
in \cite{cieliebak-frauenfelder-oancea} is based on a nonlinear deformation of
the Lagrange multiplier in Rabinowitz action functional and does not look very suitable
to compare algebraic structures on the Rabinowitz side and on the V-shaped side. In view
of this difficulty we introduce in this note an intermediate action functional,
the V-shaped action functional with delay, which shares some features of Rabinowitz action
functional and some features of the V-shaped action functional. We hope that this intermediate action functional will be able to shed some light on the ongoing scientific debate about algebraic structures in Rabinowitz Floer homology and their connection with
Tate Hochschild homology. 
\\ \\
We are considering the symplectization $\mathbb{R} \times \Sigma$ of a contact
manifold $(\Sigma,\lambda)$. We are fixing the following V-shaped function
$$h \colon (0,\infty) \to [0,\infty), \quad \rho \mapsto \rho(\ln \rho-1)+1$$
whose derivative is the logarithm so that the function $h$ attains all slopes
in a monotone increasing way. We abbreviate 
$$\mathcal{L}=C^\infty(S^1,\mathbb{R} \times \Sigma)$$
the free loop space of the symplectization where $S^1=\mathbb{R}/\mathbb{Z}$ denotes
the circle. The V-shaped action functional
with respect to the V-shaped function $h$ 
$$\mathcal{A}_0 \colon \mathcal{L} \to \mathbb{R}$$
for a loop $v=(r,x) \in \mathcal{L}$ is defined as
$$\mathcal{A}_0(v)=-\int_{S^1} v^*\lambda+\int_0^1 h(e^r)dt.$$
For the V-shaped action functional with delay 
$$\mathcal{A}_1 \colon \mathcal{L} \to \mathbb{R}$$
we simply interchange the order of
integration and applying $h$ in the second term
$$\mathcal{A}_1(v)=-\int_{S^1} v^*\lambda+h\bigg(\int_0^1 e^r dt\bigg).$$
We can interpolate between the two functionals as follows by defining for
$\theta \in [0,1]$
the functional
$$\mathcal{A}_\theta \colon \mathcal{L} \to \mathbb{R}$$
by
$$\mathcal{A}_\theta(v)=-\int_{S^1} v^* \lambda+ \theta \cdot
h\bigg(\int_0^1 e^r dt\bigg)+(1-\theta) \cdot \int_0^1 h(e^r)dt.$$
Our first lemma tells us that the critical points as well as their actions do not depend
on the parameter $\theta$ and are in natural one-to-one correspondence with
generalized periodic Reeb orbits on $\Sigma$, where the generalization corresponds to the
fact, that the period might as well be negative meaning that the orbit is traversed backward or zero meaning that the orbit is a constant point on $\Sigma$. We abbreviate by
$R$ the Reeb vector field. 
\begin{lemma}\label{l1}
Assume that $v=(r,x)$ is a critical point of $\mathcal{A}_\theta$, then 
$r$ is constant and $x$ is a solution of the ODE
\begin{equation}\label{crit}
\partial_t x(t)=rR(x(t)), \quad t \in S^1,
\end{equation}
i.e., $x$ is a periodic Reeb orbit of generalized period $r$. Moreover, the action of the critical point is given by
\begin{equation}\label{act}
\mathcal{A}_\theta(v)=1-e^r.
\end{equation}
\end{lemma}
Although the functionals have the same critical point their gradient flow lines are
different. We fix a smooth family $J_t$ for $t \in S^1$ of 
SFT-like almost complex structures and take the gradient with respect to the
$L^2$-metric obtained from the family $J_t$. Gradient flow lines of the V-shaped
action functional $\mathcal{A}_0$ are solutions $v=(r,x) \in C^\infty(
\mathbb{R} \times S^1,\mathbb{R} \times \Sigma)$ of the PDE
\begin{equation}\label{gradv}
\partial_s v(s,t)+J_t(v(s,t))\Big(\partial_t v(s,t)-r(s,t) R(v(s,t))\Big)=0.
\end{equation}
On the other hand gradient flow lines of the
V-shaped action functional with delay are solutions of the problem
\begin{equation}\label{gradvd}
\partial_s v(s,t)+J_t(v(s,t))\Bigg(\partial_t v(s,t)-\ln \bigg(\int_0^1 e^{r(s,t')}dt'\bigg) R(v(s,t))\Bigg)=0
\end{equation}
In contrast to (\ref{gradv}) this is not a PDE anymore since the integral in
front of the Reeb vector field is not local anymore but depends on the whole loop, i.e.,
this problem is a delay equation and that is the reason why we refer to $\mathcal{A}_1$
as the V-shaped action functional with delay. We rewrite the above equation equivalently
as
\begin{eqnarray*}
\partial_s v(s,t)+J_t(v(s,t))\Bigg(\partial_t v(s,t)-\tau(s) R(v(s,t))\Bigg)&=&0\\
\tau(s)=\ln \bigg(\int_0^1 e^{r(s,t')}dt'\bigg).
\end{eqnarray*}
We have the following lemma.
\begin{lemma}\label{l2}
Suppose that $v$ is a solution of the gradient flow equation (\ref{gradvd})
converging asymptotically to critical points $v_\pm=\lim_{s \to \pm \infty} v(s)$.
Then for every $s \in \mathbb{R}$ it holds that 
$$\partial_s \tau(s) \geq -\frac{\mathcal{A}_1(v_-)-\mathcal{A}_1(v_+)}{1-\mathcal{A}_1(v_+)}>-1.$$
\end{lemma}
It follows from Lemma~\ref{l2} by standard arguments, see for instance \cite{cieliebak-frauenfelder-oancea}, that if $v=(r,x)$ is as in the Lemma, then
\begin{equation}\label{lapl}
\Delta r\geq -\frac{\mathcal{A}_1(v_-)-\mathcal{A}_1(v_+)}{1-\mathcal{A}_1(v_+)} > -1,
\end{equation}
i.e., the Laplacian of $r$ is uniformly bounded from below. This guarantees an upper
bound of $r$ just in terms of asymptotic datas of $r$.
\\ \\
We as well recall Rabinowitz action functional on a symplectization
$$\mathcal{A}_2 \colon \mathcal{L} \times \mathbb{R} \to \mathbb{R}$$
which for $(v,\tau)=(r,x,\tau) \in \mathcal{L} \times \mathbb{R}$
is given by
$$\mathcal{A}_2(v,\tau)=-\int_{S^1} v^* \lambda+\tau \int_0^1 e^r dt-\tau.$$ 
If one interprets $\tau$ as a Lagrange multiplier, than $\mathcal{A}_2$
is the Lagrange multiplier functional of the negative area functional
with respect to the constraint given by the mean value of the function
$e^r-1$. By the theory of Lagrange multipliers critical points of $\mathcal{A}_2$
correspond to critical points of the restriction of the negative area functional
to
$$\overline{\mathcal{L}}=\bigg\{v=(r,x) \in \mathcal{L}: \int_0^1 e^r dt=1\bigg\},$$
namely
$$\mathcal{A}_3 \colon \overline{\mathcal{L}} \to \mathbb{R}, \quad
v \mapsto -\int_{S^1}v^*\lambda.$$
This functional was discussed by the author in \cite{frauenfelder}.
\\ \\
The gradient flow equations of
the Rabinowitz action functionals $\mathcal{A}_2$ and $\mathcal{A}_3$ are delay
equations as well. For $\mathcal{A}_2$ they are solutions 
$(v,\tau)=(r,x,\tau) \in C^\infty(\mathbb{R}\times S^1,\mathbb{R} \times \Sigma)
\times C^\infty(\mathbb{R},\mathbb{R})$ of the problem
\begin{eqnarray}\label{gradrab1}
\partial_s v(s,t)+J_t(v(s,t))\Big(\partial_t v(s,t)-\tau(s)R(v(s,t))\Big)&=&0\\ \nonumber
\partial_s \tau(s)+\int_0^1 e^{r(s,t')}dt'&=&1,
\end{eqnarray}
and for $\mathcal{A}_3$ solutions $v=(r,x) \in C^\infty(\mathbb{R} \times S^1,
\mathbb{R} \times \Sigma)$ of the problem
\begin{eqnarray}\label{gradrab2}
\partial_s v(s,t)+J_t(v(s,t))\Big(\partial_t v(s,t)+\mathcal{A}_3(v_s)R(v(s,t))\Big)&=&0\\ \nonumber
\int_0^1 e^{r(s,t')}dt'&=&1,
\end{eqnarray}
where we abbreviate by $v_s \in \overline{\mathcal{L}}$ the loop $t \mapsto v(s,t)$, 
see \cite{frauenfelder}.
\\ \\
All four problems have the form
\begin{equation}\label{percr}
\partial_s v+J(v)\big(\partial_t v-\tau R(v)\big)=0,
\end{equation}
i.e., are perturbed Cauchy-Riemann equations where the perturbation is in the Reeb direction. In the last three problems, $\tau$ only depends on $s$ and is independent of
$t$, while in the first problem it depends on both variables. We see that the $t$-independence of $\tau$ is a feature the gradient
flow equation for the V-shaped action functional with delay has in common with the
gradient flow equations of the Rabinowitz action functionals. On the other hand
we have the following ascending chain of spaces
$$\overline{\mathcal{L}} \subset \mathcal{L} \subset \mathcal{L} \times \mathbb{R}$$
where the codimension in each step is one. Both V-shaped action functionals live on the
intermediate space where the Rabinowitz action functionals live on the two extremal 
spaces. 
\\ \\
We have a natural projection
$$\Pi \colon \mathcal{L} \to \overline{\mathcal{L}}, \quad
(r,x) \mapsto \bigg(r-\ln \int_0^1 e^r dt,x\bigg). $$
In case of $\mathcal{A}_3$, i.e., the restriction of the negative area functional
to the constraint $\overline{\mathcal{L}}$, the Lagrange multiplier $\tau(s)$ is uniquely
determined by the loop $v_s \in \overline{\mathcal{L}}$.
Therefore solutions of problem (\ref{gradrab2}) can be uniquely characterized as solutions
of the perturbed Cauchy-Riemann equation (\ref{percr}) where $\tau$ only depends on
$s$ but not on $t$ and $v_s \in \overline{\mathcal{L}}$ for every $s \in \mathbb{R}$. 
We abbreviate by $\mathcal{M}_i$ for $0 \leq i \leq 3$ the moduli spaces of finite
energy solutions of the gradient flow lines of the four functionals
$\mathcal{A}_i$. We have natural maps
$$\Pi_* \colon \mathcal{M}_1 \to \mathcal{M}_3, \,\,\Pi_* \colon
\mathcal{M}_2 \to \mathcal{M}_3, \quad v \mapsto \Pi v.$$
This is not true for $\mathcal{M}_0$ since in this case $\tau$ usually depends on
$t$ as well. In \cite{frauenfelder} it was shown that 
$\Pi_* \colon \mathcal{M}_2 \to \mathcal{M}_3$ is a bijection, i.e.,
there is a natural bijection between finite energy gradient flow lines
of the two Rabinowitz action functionals. The following theorem
tells us that the same is true for the gradient flow lines of
the V-shaped action functional with delay
\begin{thm}\label{l3}
The map $\Pi_* \colon \mathcal{M}_1 \to \mathcal{M}_3$ is a bijection.
\end{thm}
In particular, we see from the lemma combined with \cite{frauenfelder} that there is a natural bijection between
the finite energy gradient flow lines of the V-shaped action functional with delay
and the ones of the Rabinowitz action functional $\mathcal{A}_2$. 
\\ \\
\emph{Acknowledgements: } The author acknowledges partial support by DFG grant
FR 2637/2-2.

\section{Proof of Lemma~\ref{l1} and derivation of the gradient flow equation}

Before we embark on the proof of Lemma~\ref{l1} we recall some notation. 
We assume that $\Sigma=\Sigma^{2n-1}$ is a closed odd dimensional manifold
of dimension $2n-1$ and $\lambda \in \Omega^1(\Sigma)$ is a contact form on
$\Sigma$, i.e., 
$$\lambda \wedge (d\lambda)^{n-1}>0$$ 
is a volume form on $\Sigma$. We abbreviate by $\xi$ the contact structure, i.e.,
the hyperplane distribution
$$\xi=\ker \lambda$$
and by $R$ the Reeb vector field on $\Sigma$ implicitly defined by the conditions
$$\lambda(R)=1, \quad d\lambda(R,\cdot)=0.$$
By abuse of notation we use the same letter $\lambda$ for the extension
of the one-form to the symplectization
$\mathbb{R} \times \Sigma$ of $\Sigma$, namely we set
$$\lambda_{r,x}=e^r \lambda_x, \quad (r,x) \in \mathbb{R} \times \Sigma.$$
The exterior derivative of $\lambda$
$$\omega_{r,x}=d\lambda_{r,x}=e^r dr \wedge \lambda_x+e^r d\lambda_x$$
is then a symplectic form on $\mathbb{R} \times \Sigma$. The bundle $\xi$ and
the Reeb vector field $R$ we extend trivially to the symplectization and use there
by abuse of notation the same letters again. 
\\ \\
\textbf{Proof of Lemma~\ref{l1}: } Suppose that $v=(r,x) \in \mathcal{L}$ and
$\widehat{v}$ is a vector field along $v$, i.e., a tangent vector
$$\widehat{v}=(\widehat{r},\widehat{x}) \in T_v \mathcal{L}=
\Gamma\big(v^*T(\mathbb{R} \times \Sigma)\big)=C^\infty(S^1,\mathbb{R})
\times \Gamma(x^*T\Sigma).$$
In the following computation we denote by $L_{\widehat{v}}$ the Lie derivative
in direction of $\widehat{v}$ and use Cartan's formula $L_{\widehat{v}}
=d\iota_{\widehat{v}}+\iota_{\widehat{v}}d$.
\begin{eqnarray*}
d\mathcal{A}_\theta(v)\widehat{v}&=&-\int_{S^1} v^*L_{\widehat{v}} \lambda
+\theta \cdot h'\bigg(\int_0^1 e^{r(t)}dt\bigg)\int_0^1 e^{r(t)} \widehat{r}(t)dt\\
& &+(1-\theta)\cdot \int_0^1 h'\big(e^{r(t)}\big)e^{r(t)} \widehat{r}(t)dt\\
&=&-\int_{S^1} v^* d \iota_{\widehat{v}} \lambda-\int_{S^1} v^* \iota_{\widehat{v}}
d \lambda\\
& &+\theta \cdot \ln \bigg(\int_0^1 e^{r(t)}dt\bigg)\int_0^1
e^{r(t)} \lambda_{x(t)} \big(R(x(t))\big) \widehat{r}(t)dt\\
& &+(1-\theta)\cdot \int_0^1 \ln \big(e^{r(t)}\big) e^{r(t)}
\lambda_{x(t)} \big(R(x(t))\big) \widehat{r}(t)dt\\
&=&-\int_{S^1} d v^* \iota_{\widehat{v}} \lambda-\int_{S^1} v^* \iota_{\widehat{v}}
\omega+\theta \cdot \ln \bigg(\int_0^1 e^{r}dt\bigg)
\int_0^1 \omega\big(\widehat{v},R(v)\big)dt\\
& &+(1-\theta) \cdot \int_0^1 r\omega\big(\widehat{v},R(v)\big)dt\\
&=&-\int_0^1 \omega(\widehat{v},\partial_t v)dt+\int_0^1 \omega
\Bigg(\widehat{v},\theta\cdot \ln\bigg(\int_0^1 e^{r}dt\bigg)R(v)\Bigg)dt\\
& &+\int_0^1 \omega \Big(\widehat{v},(1-\theta)\cdot r R(v)\Big)dt\\
&=&\int_0^1 \omega\Bigg(\partial_t v-\theta \cdot \ln \bigg(\int_0^1 e^r dt\bigg)R(v)
-(1-\theta)\cdot r R(v), \widehat{v}\Bigg)dt.
\end{eqnarray*}
For later reference we summarize this to
\begin{equation}\label{diff}
d\mathcal{A}_\theta(v)\widehat{v}=\int_0^1 \omega\Bigg(\partial_t v-\bigg(\theta \cdot \ln \bigg(\int_0^1 e^r dt\bigg)
+(1-\theta)\cdot r \bigg)R(v),\widehat{v}\Bigg)dt.
\end{equation} 
In particular, we infer from (\ref{diff}) that $v \in \mathrm{crit}(\mathcal{A}_\theta)$
if and only if $v$ is a solution of the problem
\begin{equation}\label{crit0}
\partial_t v=\bigg(\theta \cdot \ln \bigg(\int_0^1 e^r dt\bigg)
+(1-\theta)\cdot r \bigg)R(v).
\end{equation}
Since $dr(R)=0$ we infer from (\ref{crit0}) that $r$ is constant. Therefore
the scaling factor in front of the Reeb vector field simplifies as follows
$$\theta \cdot \ln \bigg(\int_0^1 e^r dt\bigg)
-(1-\theta)\cdot r=\theta \cdot \ln e^r 
+(1-\theta)\cdot r =\theta \cdot r+(1-\theta)\cdot r=r.$$
This shows that problem (\ref{crit0}) is equivalent to (\ref{crit}) and the first
assertion of the lemma is proved.
\\ \\
It remains to prove the second assertion, namely the formula (\ref{act}) for 
the action of a critical point. Plugging (\ref{crit}) into the definition of
$\mathcal{A}_\theta$ we compute using that $r$ is independent of time
\begin{eqnarray*}
\mathcal{A}_\theta(v)&=&-\int_0^1 e^r \lambda_x\big(rR) dt+\theta h\big(e^r\big)+(1-\theta)
h\big(e^r\big)\\
&=&-re^r+h\big(e^r)\\
&=&-re^r+e^r(r-1)+1\\
&=&1-e^r.
\end{eqnarray*}
This shows (\ref{act}) and hence the lemma is proved. \hfill $\square$
\\ \\
We derive the gradient flow equation for the functionals $\mathcal{A}_\theta$
simultaneously for every $\theta \in [0,1]$. The gradient flow equations
(\ref{gradv}) and (\ref{gradvd}) follow than by specializing to the cases
$\theta=0$, respectively $\theta=1$.  To write down the gradient flow
equation we have to choose a metric on the free loop space. For that purpose we
choose a smooth family of SFT-like almost complex structures $J_t$ and take the
$L^2$-metric with respect to this family. SFT-like almost complex structures
are used in Symplectic Field Theory \cite{eliashberg-givental-hofer} and we recall their definition. The contact condition implies that the vector
bundle $(\xi,d\lambda|_\xi) \to \Sigma$ is symplectic. Hence we choose a family
of $d\lambda|_\xi$ compatible almost complex structures $J_t$ on $\xi$, i.e.,
$d\lambda|_\xi(\cdot, J_t \cdot)$ is a bundle metric on $\xi$. 
We extend this family canonically to a family of $\mathbb{R}$-invariant
almost complex structures on the tangent space of the symplectization
by requiring that $J_t$ interchanges the Reeb vector field $R$ and the Liouville
vector field $\partial_r$, namely 
$$J_tR=-\partial_r, \quad J_t\partial_r=R.$$
In particular, $J_t$ is compatible with the symplectic form $\omega$ in the sense that
$\omega(\cdot,J_t \cdot)$ is a Riemannian metric on $\mathbb{R} \times \Sigma$. 
Following Floer \cite{floer} we use this family of $\omega$-compatible almost
complex structures to define an $L^2$-metric $g=g_J$ on the free loop space $\mathcal{L}$. If
$v \in \mathcal{L}$ and $\widehat{v}_1, \widehat{v}_2 \in T_v \mathcal{L}$ are
tangent vectors at $v$, i.e., vector fields along $v$, we define the $L^2$-inner product of $\widehat{v}_1$ and $\widehat{v}_2$ by
$$g\big(\widehat{v}_1,\widehat{v}_2\big)=\int_0^1 \omega\Big(\widehat{v}_1(t),
J_t\big(v(t)\big)\widehat{v}_2(t)\Big)dt.$$
We define the $L^2$-gradient $\nabla \mathcal{A}_\theta(v)$ 
of $\mathcal{A}_\theta$ at $v \in \mathcal{L}$
implicitly by the requirement that
$$d \mathcal{A}_\theta(v) \widehat{v}=g\big(\nabla \mathcal{A}_\theta(v),\widehat{v}\big),
\quad \forall\,\,\widehat{v} \in T_v \mathcal{L}.$$
We claim that
$$\nabla \mathcal{A}_\theta(v)(t)=J_t(v(t))\Bigg(\partial_t v-\bigg(\theta \cdot \ln \bigg(\int_0^1 e^r dt\bigg)
+(1-\theta)\cdot r \bigg)R(v)\Bigg)(t), \quad t \in S^1.$$
To see that we rewrite (\ref{diff}) using that $J^2=-\mathrm{id}$
\begin{eqnarray*}
d \mathcal{A}_\theta(v) \widehat{v}&=&
\int_0^1 \omega\Bigg(\widehat{v}, J(v)^2\Bigg(\partial_t v-\bigg(\theta \cdot \ln \bigg(\int_0^1 e^r dt\bigg)
+(1-\theta)\cdot r \bigg)R(v)\Bigg)\Bigg)dt\\
&=&g\Bigg(\widehat{v}, J(v)\Bigg(\partial_t v-\bigg(\theta \cdot \ln \bigg(\int_0^1 e^r dt\bigg)
+(1-\theta)\cdot r \bigg)R(v)\Bigg)\Bigg),
\end{eqnarray*}
which implies the formula above for the gradient. A gradient flow line
is formally a smooth map $v \in C^\infty(\mathbb{R},\mathcal{L})$ solving
the ``ODE" on $\mathcal{L}$
$$\partial_s v(s)+\nabla \mathcal{A}_\theta(v)(s)=0, \quad s \in \mathbb{R},$$
which we interpret as a smooth map
$v \in C^\infty(\mathbb{R} \times S^1, \mathbb{R} \times \Sigma)$ solving the
problem
\begin{equation}\label{grad}
\partial_s v+J(v)\Bigg(\partial_t v-\bigg(\theta \cdot \ln \bigg(\int_0^1 e^r dt\bigg)
+(1-\theta)\cdot r \bigg)R(v)\Bigg)=0.
\end{equation}
Let us spell out in detail the extremal cases $\theta=0$ and $\theta=1$. For
$\theta=0$ problem (\ref{grad}) becomes
\begin{equation}\label{grad0}
\partial_s v(s,t)+J_t(v(s,t))\Big(\partial_t v(s,t)-r(s,t)R(v(s,t))\Big)=0, \quad
(s,t) \in \mathbb{R} \times S^1.
\end{equation}
This is a PDE on the cylinder. In contrast to (\ref{grad0}) the problem (\ref{grad})
is not a PDE anymore for positive $\theta$, since in this case the factor in front
of the Reeb vector field not just depends on $s$ and $t$ but on the whole loop
$$r_s \colon S^1 \to \mathbb{R}, \quad t' \mapsto r(s,t').$$
For $\theta=1$ this becomes for $(s,t) \in \mathbb{R} \times S^1$
\begin{equation}\label{grad1}
\partial_s v(s,t)+J_t(v(s,t))\Bigg(\partial_t v(s,t)-\ln\bigg(\int_0^1 e^{r(s,t')}dt'\bigg)
R(v(s,t))\Bigg)=0. 
\end{equation}
For every $\theta \in [0,1]$ the problem (\ref{grad}) is a perturbation of
the Cauchy-Riemann equation in the direction of the Reeb vector field. 
In case of (\ref{grad1}) the scaling factor in front of the Reeb vector field only
depends on $s$ and not on $t$. We abbreviate this scaling factor for
$(s,t) \in \mathbb{R} \times S^1$ by
$$\tau(s,t)=\theta \cdot \ln \bigg(\int_0^1 e^{r(s,t')}dt'\bigg)
+(1-\theta) \cdot r(s,t).$$
With this notion we can write (\ref{grad}) more compactly as
$$\partial_s v+J(v)\big(\partial_t v-\tau R(v)\big)=0.$$

\section{Proof of Lemma~\ref{l2}}

Suppose that $v$ is a gradient flow line of the V-shaped action functional
with delay $\mathcal{A}_1$
$$\partial_s v+\nabla \mathcal{A}_1(v)=0,$$
i.e., $v=(r,x) \in C^\infty(\mathbb{R} \times S^1, \mathbb{R} \times \Sigma)$
is a solution of (\ref{gradvd}). We abbreviate
$$T(s):=\int_0^1 e^{r(s,t)}dt$$
and 
$$\tau(s):=\ln(T(s))=\ln\bigg(\int_0^1 e^{r(s,t)}dt\bigg).$$
Since $J_t$ is SFT-like we obtain from (\ref{gradvd})
$$\partial_s r-\lambda_x(\partial_t x)+\tau=0.$$
Abbreviating further $v_s$ the loops obtained by evaluating $v$ at time $s$ we compute
\begin{eqnarray*}
\partial_s T(s)&=&\int_0^1 e^{r(s,t)} \partial_s r(s,t)dt\\
&=&\int_0^1 e^{r(s,t)} \Big(\lambda_{x(s)}\big(\partial_t x(s,t)\big)-\tau(s)\Big)dt\\
&=&\int_{S^1} v_s^*\lambda-\tau(s) T(s)\\
&=&-\mathcal{A}_1(v_s)+h\bigg(\int_0^1 e^{r(s,t)}dt\bigg)-\ln(T(s)) T(s)\\
&=&-\mathcal{A}_1(v_s)+h(T(s))-\ln(T(s)) T(s)\\
&=&-\mathcal{A}_1(v_s)+T(s) \ln(T(s))-T(s)+1-\ln(T(s))T(s)\\
&=&-\mathcal{A}_1(v_s)-T(s)+1
\end{eqnarray*}
and therefore
\begin{equation}\label{der}
\partial_s \tau(s)=\frac{1-\mathcal{A}_1(v_s)}{e^{\tau(s)}}-1
\end{equation}
In the case where the gradient flow line is constant the assertion of the lemma
is trivial. Hence we assume in the following that $v$ is nonconstant. 
Since $v$ is a gradient flow line $\mathcal{A}_1$ it follows that
$\mathcal{A}_1$ is strictly decreasing along $v$ so that we have
$$\mathcal{A}_1(v_+)<\mathcal{A}_1(v_s)<\mathcal{A}_1(v_-)<1.$$
The last inequality follows from Lemma~\ref{l1}. We rewrite this equivalently as
\begin{equation}\label{acin}
1-\mathcal{A}_1(v_+)>1-\mathcal{A}_1(v_s)>1-\mathcal{A}_1(v_-)>0
\end{equation}
We claim that
\begin{equation}\label{T}
e^{\tau(s)} \leq 1-\mathcal{A}_1(v_+).
\end{equation}
We first show (\ref{T}) under the assumption that $\tau$ assumes its global minimum,
i.e., there exists $s_0 \in \mathbb{R}$ such that
$$\tau(s_0)=\min \tau.$$
In particular, we have
$$\partial_s \tau(s_0)=0$$
and therefore in view of (\ref{der})
$$e^{\tau(s_0)}=1-\mathcal{A}_1(v_{s_0})$$
implying that
$$\tau(s_0)=\ln\big(1-\mathcal{A}(v_{s_0})\big)\geq \ln\big(1-\mathcal{A}(v_+)\big).$$
Since $\tau$ assumed its minimum at $s_0$ we have in particular for every $s \in \mathbb{R}$
$$\tau(s) \geq \ln\big(1-\mathcal{A}(v_+)\big)$$
and therefore
$$e^{\tau(s)} \leq 1-\mathcal{A}(v_+).$$
This proves (\ref{T}) in the case where $\tau$ assumes its global minimum. 
\\ \\
It remains to discuss the case where $\tau$ does not attend its minimum. We first
show that in this case the derivative of $\tau$ never vanishes. Suppose by contradiction
that there exists $s_0 \in \mathbb{R}$ such that
$$\partial_s \tau(s_0)=0.$$
We infer from (\ref{der}) that 
$$\partial^2 \tau(s_0)=-\frac{\partial_s \mathcal{A}_1(v_{s_0})}{e^{\tau(s_0)}}> 0$$
where the inequality follows from the fact that $\mathcal{A}_1$ is strictly decreasing
along its nonconstant gradient flow line. This implies that $s_0$ is a local minimum
of $\tau$ and $\tau$ cannot have local maxima. In particular $s_0$ has to be the global minimum in contradiction to our assumption that this does not exist. It follows that
$$\partial_s\tau(s) \neq 0,\quad \forall\,\,s \in \mathbb{R},$$
i.e.,\,$\tau$ is strictly monotone. Abbreviating $\tau_\pm =\lim_{s \to \pm \infty}
\tau(s)$ we must have
$$\tau_+=\ln\big(1-\mathcal{A}_1(v_+)\big) >\tau(s)> \ln \big(1-\mathcal{A}_1(v_-)\big)
=\tau_-$$
and therefore
$$e^{\tau(s)} \leq 1-\mathcal{A}_1(v_+).$$
This proves that (\ref{T}) holds as well in the case where $\tau$ does not attain
its global minimum and the formula is proved. 
\\ \\
Since $1-\mathcal{A}(v_s)>0$ by (\ref{acin}) we obtain from (\ref{der}) and (\ref{T}) 
$$\partial_s \tau(s) \geq\frac{1-\mathcal{A}(v_s)}{1-\mathcal{A}(v_+)}-1
= -\frac{\mathcal{A}(v_s)-\mathcal{A}(v_+)}{1-\mathcal{A}(v_+)}
\geq -\frac{\mathcal{A}(v_-)-\mathcal{A}(v_+)}{1-\mathcal{A}(v_+)}>-1.$$
This finishes the proof of the lemma. 

\section{Proof of Theorem~\ref{l3}}

We construct a map $\mathcal{R} \colon \mathcal{M}_3 \to \mathcal{M}_1$ inverse
to $\Pi_* \colon \mathcal{M}_1 \to \mathcal{M}_3$. Suppose that $v=(r,x) \in \mathcal{M}_3$, i.e., a finite energy solution of (\ref{gradrab2}). We are looking for
$\rho \in C^\infty(\mathbb{R},\mathbb{R})$ such that
$$w=\rho_* v=(r+\rho,x)$$
is a solution of (\ref{gradvd}), i.e., a gradient flow line of $\mathcal{A}_1$.
Since $J_t$ is SFT-like and in particular $R$-invariant, $w$ is a solution of the
problem
$$\partial_s w+J(w)\Big(\partial_t w +\big(\mathcal{A}_3(v)+\partial_s \rho\big)R(w)\Big)=0.$$
In order that this becomes a solution of (\ref{gradvd}) we need that for
every $s \in \mathbb{R}$ it holds that
\begin{eqnarray*}
\mathcal{A}_3(v_s)+\partial_s \rho(s)&=&\ln\bigg(\int_0^1 e^{r(s,t)+\rho(s)}dt\bigg)\\
&=&\ln\bigg(e^{\rho(s)}\int_0^1 e^{r(s,t)}dt\bigg)\\
&=&\rho(s)+\ln\bigg(\int_0^1 e^{r(s,t)}dt\bigg)\\
&=&\rho(s),
\end{eqnarray*}
where in the last equality we have used that $v_s \in \overline{\mathcal{L}}$, since
$v \in \mathcal{M}_3$. We see that $\rho$ has to be a solution of the ODE
\begin{equation}\label{req}
\partial_s \rho(s)=\rho(s)-\mathcal{A}_3(v_s).
\end{equation}
Abbreviate 
$$\sigma(s)=\rho(s)-\mathcal{A}_3(v_s).$$
From (\ref{req}) we obtain
$$\partial_s \sigma(s)=\partial_s \rho(s)-\partial_s \mathcal{A}_3(v_s)
=\rho(s)-\mathcal{A}_3(v_s)-\partial_s \mathcal{A}_3(v_s)=\sigma(s)-\partial
_s \mathcal{A}_3(v_s)$$
so that $\sigma$ is a solution of the ODE
\begin{equation}\label{seq}
\partial_s \sigma(s)=\sigma(s)-\partial_s \mathcal{A}_3(v_s).
\end{equation}
Since $v$ is a gradient flow line of $\mathcal{A}_3$ we have
$\partial_s \mathcal{A}_3(v) \geq 0$, and since the gradient flow line has
finite energy the integral of $\partial_s \mathcal{A}_3(v)$ is finite. In particular,
$$\partial_s \mathcal{A}_3(v) \in L^1(\mathbb{R}).$$
Since it is smooth as well we also have
$$\partial_s \mathcal{A}_3(v) \in L^2(\mathbb{R}).$$ 
The operator
$$D \colon W^{1,2}(\mathbb{R}) \to L^2(\mathbb{R}), \quad
\sigma \mapsto \partial_s \sigma-\sigma$$
is an injective Fredholm operator of index zero, i.e., an isomorphism between
the two Hilbert spaces. Therefore there exists a unique
$\sigma \in W^{1,2}(\mathbb{R})$ solving (\ref{seq}), namely
$$\sigma=D^{-1}\big(-\partial_s \mathcal{A}_3(v)\big).$$
Now set
$$\rho_v:=\mathcal{A}_3(v)+D^{-1}\big(-\partial_s \mathcal{A}_3(v)\big)$$
then $w:=(\rho_v)_*(v)$ is a finite energy solution of (\ref{gradvd}), i.e., an
element of the moduli space $\mathcal{M}_1$, and we have a well-defined map
$$\mathcal{R} \colon \mathcal{M}_3 \to \mathcal{M}_1, \quad v \mapsto (\rho_v)_* v.$$
It remains to check that $\mathcal{R}$ is inverse to $\Pi_*$. We carry this out in two steps.
\\ \\
\textbf{Step\,1: } We have $\Pi_* \circ \mathcal{R}=\mathrm{id} \colon
\mathcal{M}_3 \to \mathcal{M}_3$, i.e., $\mathcal{R}$ is right inverse
to $\Pi_*$. 
\\ \\
This is the same argument as the proof of Step\,1 of Theorem\,5.1 in \cite{frauenfelder}. 
Assume that $v=(r,x) \in \mathcal{M}_3$. Since $\Pi_*$ and $\mathcal{R}$
are both $t$-independent translations in the $\mathbb{R}$-direction there exists
$\rho \in C^\infty(\mathbb{R},\mathbb{R})$ such that
$$w:=\Pi_* \circ \mathcal{R}(v)=\rho_* v=(r+\rho,x).$$
Since both $v$ and $w$ belong to $\mathcal{M}_3$ we must have for every $s \in 
\mathbb{R}$
$$\int_0^1 e^{r(s,t)}dt=1, \quad \int_0^1 e^{r(s,t)+\rho(s)}dt=1$$
implying that
$$\rho=0.$$
Therefore $w=v$ so that we have
$$\Pi_* \circ \mathcal{R}(v)=v.$$
Since $v \in \mathcal{M}_3$ was arbitrary Step\,1 follows. 
\\ \\
\textbf{Step\,2: } We have $\mathcal{R} \circ \Pi_*=\mathrm{id} \colon
\mathcal{M}_1 \to \mathcal{M}_1$, i.e., $\mathcal{R}$ is left inverse
to $\Pi_*$. 
\\ \\
Assume that $v=(r,x) \in \mathcal{M}_1$. Again since $\mathcal{R}$ and $\Pi_*$
are both $t$-independent translations in the $\mathbb{R}$-direction there exists
$\rho \in C^\infty(\mathbb{R},\mathbb{R})$ such that
$$w:=\mathcal{R} \circ \Pi_*(v)=\rho_* v=(r+\rho,x).$$
Since $v$ solves (\ref{gradvd}) we have
\begin{eqnarray*}
0&=&\partial_s w+J(w)\Bigg(\partial_t w+\bigg(\partial_s \rho-\ln \int_0^1 e^{r}dt'\bigg)
R(w)\Bigg)\\
&=&\partial_s w+J(w)\Bigg(\partial_t w+\bigg(\partial_s \rho+\rho-\ln \int_0^1 e^{r+\rho}dt'\bigg)R(w)\Bigg).
\end{eqnarray*}
Since $w$ is a solution of problem (\ref{gradvd}) as well we must have
$$\partial_s \rho+\rho=0$$
i.e.,
$$\rho(s)=\rho_0e^{-s}.$$
Since $w$ has finite energy it follows that $\rho_0=0$ and therefore
$$\rho=0.$$
This proves that $w=v$ and hence
$$\mathcal{R} \circ \Pi_*(v)=v.$$
Since $v \in \mathcal{M}_3$ was arbitrary Step\,2 follows and the Theorem is proved.

\end{document}